\documentclass{amsart}
\usepackage{amssymb}\usepackage{amsmath}\usepackage{a4}
\def\operatorname#1{\hbox{#1}}\def\text#1{\hbox{#1}}
\newtheorem{theorem}{Theorem}[section]
\newtheorem{lemma}[theorem]{Lemma}

\newtheorem{definition}[theorem]{Definition}
\def\qedbox{\hbox{$\rlap{$\sqcap$}\sqcup$}}
\makeatletter
  
  \@addtoreset{equation}{section}
 \makeatother
\begin{document}
\title{Pseudo-Riemannian Jacobi--Videv Manifolds}

\author{P. Gilkey and S. Nik\v cevi\'c}
\address{PG: Mathematics Department, University of Oregon,
 Eugene Oregon 97403 USA.}
\email{gilkey@uoregon.edu}
\address{SN: Mathematical Institute, SANU,
 Knez Mihailova 35, p.p. 367,
 11001 Belgrade,
 Serbia.}
\email{stanan@mi.sanu.ac.yu}

\begin{abstract} We exhibit several families of Jacobi--Videv pseudo-Riemannian manifolds which are not
Einstein. We also exhibit Jacobi--Videv algebraic curvature tensors where the Ricci operator defines an almost complex structure.
\end{abstract}

\maketitle
\def\BB{\mathcal{B}}
\def\MM{\mathfrak{M}}
\def\Gr{\operatorname{Gr}}
\def\qedbox{\hbox{$\rlap{$\sqcap$}\sqcup$}}
\def\bork{\vphantom{\vrule height12pt}}
\section{Introduction}

Studying additional algebraic properties of the curvature tensor and relating these properties to the underlying geometry is an active field of
investigation recently. Although Lorentzian geometry plays a central role in mathematical physics, the higher
signature context is important as well. We refer to \cite{CGM05, M05} for a discussion of signature $(2,2)$ Walker metrics; these are manifolds which
admit a parallel totally isotropic $2$-plane field. Dual and anti-self dual metrics are discussed by \cite{BS, HSW} in the higher
signature setting. Manifolds with a nilpotent Ricci operator of higher order appear naturally \cite{DM}, and manifolds of signatures other than
Riemannian or Lorentzian are important in Brane theory
\cite{LPS}.

In this paper, we shall study when the Ricci operator and the Jacobi operator commute. Let $\nabla$ be the Levi-Civita connection of a
pseudo-Riemannian manifold $(M,g)$ of signature $(p,q)$; we shall primarily be concerned with the case $p\ge1$ and $q\ge1$. Let
$\mathcal{R}(x,y):=\nabla_x\nabla_y-\nabla_y\nabla_x-\nabla_{[x,y]}$ be the Riemann curvature operator, let
$\mathcal{J}(x):y\rightarrow\mathcal{R}(y,x)x$ be the Jacobi operator, and let $\rho$ be the Ricci operator. Following seminal work of
Videv, we say that
$\mathcal{M}$ is {\it Jacobi--Videv} if $\mathcal{J}(x)\rho=\rho\mathcal{J}(x)$ for all $x$.

Clearly if $\mathcal{M}$ is Einstein, i.e. if $\rho=c\operatorname{id}$, then $\mathcal{M}$ is Jacobi--Videv. If $\mathcal{M}$ is
indecomposable, the converse implication holds in the Riemannian setting; any indecomposable Riemannian Jacobi--Videv manifold is necessarily
Einstein \cite{GP06}; we also refer to related work \cite{IZV06}. This implication fails in the indefinite context.
One has the following family of examples which are Jacobi--Videv and not Einstein. Manifolds in this family have been studied previously in different
contexts, see for example
\cite{BD00, BV97, CLPT90, CGV05,GN05}; we also refer to \cite{BGGV07a,SV92}

\begin{definition}\label{defn-1.1}
\rm Let $k\ge1$, let $\ell\ge1$, and $m=2k+\ell$. Introduce coordinates 
$$(x_1,...,x_k,y_1,...,y_\ell,\bar x_1,...,\bar x_k)\quad\text{on}\quad\mathbb{R}^m\,.$$
Let indices $i,j$ range from $1$ through $k$ and index the collections $\{\partial_{x_1},...,\partial_{x_k}\}$ and
$\{\partial_{\bar x_1},...,\partial_{\bar x_k}\}$. Let indices $a,b$ range from $1$ through $\ell$ and index the collection
$\{\partial_{y_1},...,\partial_{y_\ell}\}$. Let $S^2(\mathbb{R}^k)$ be the set of symmetric $2$-tensors on $\mathbb{R}^k$ and let $\psi$ be
a smooth map from $\mathbb{R}^{\ell}$ to $S^2(\mathbb{R}^k)$. Let $C_{ab}=C_{ba}$ define a non-degenerate inner product
of signature
$(u,v)$ on $\mathbb{R}^\ell$ where $u+v=\ell$. We use
$\psi$ and
$C$ to define a pseudo-Riemannian manifold $\mathcal{M}=\mathcal{M}_{C,\psi}:=(\mathbb{R}^{2k+\ell},g_{C,\psi})$ where $g=g_{C,\psi}$ is the
pseudo-Riemannian manifold of signature
$(k+u,k+v)$ with non-zero components
$$g(\partial_{x_i},\partial_{x_j}):=-2\psi_{ij}(\vec y),\quad
 g(\partial_{y_a},\partial_{y_b})=C_{ab},\quad
 g(\partial_{x_i},\partial_{\bar x_i})=1\,.
$$
\end{definition}

\begin{theorem}\label{thm-1.2}
The manifold $\mathcal{M}_{C,\psi}$ of Definition \ref{defn-1.1} is Jacobi--Videv with nilpotent Ricci operator $\rho$. Furthermore,
$\mathcal{M}_{C,\psi}$ is Einstein if and only if for any $i,j$ with $1\le i,j\le k$ we have
$\sum_{ab}C^{ab}\partial_{y_a}\partial_{y_b}\psi_{ij}=0$.
\end{theorem}

We note that if $\psi$ is a periodic function, then the metric $g$ descends to
define a metric on the torus $\mathbb{T}^{2k+\ell}$. Thus there are compact examples of Jacobi--Videv manifolds which are not Einstein.

One says that a pseudo-Riemannian manifold $\mathcal{M}$ is {\it curvature homogeneous} if given any two points $P$ and $Q$ of $M$,
there is an isometry $\phi_{P,Q}$ from $T_PM$ to $T_QM$ so that $\phi_{P,Q}^*R_Q=R_P$.

Although a pseudo-Riemannian manifold need not be Einstein, it is known \cite{GP06} that if $\mathcal{M}$ is
an indecomposable Jacobi--Videv manifold, then either $\rho$ has only one real eigenvalue
or $\rho$ has two complex eigenvalues which are complex conjugates; such a manifold is said to be {\it pseudo-Einstein}. Clearly if $\rho$ is
nilpotent, then
$0$ is the only eigenvalue of
$\rho$. This does not, however, imply $\mathcal{M}$ is Jacobi--Videv as the following result shows:

\begin{theorem}\label{thm-1.3}
Let $\{x,y,z,\bar x\}$ be coordinates on $\mathbb{R}^4$. Let $\phi\in C^\infty(\mathbb{R})$. Assume $\phi^{\prime\prime}(y)\ne0$ for all $y$. Let
$\mathcal{M}:=(\mathbb{R}^4,g)$ where
$
g(\partial_{x},\partial_{\bar x})=g(\partial_y,\partial_y)=g(\partial_z,\partial_z)=1$ and 
$g(\partial_{x},\partial_z)=2\phi(y)$. Then:
\begin{enumerate}
\item $\operatorname{Rank}\{\rho\}=3$, $\operatorname{Rank}\{\rho^2\}=2$, $\operatorname{Rank}\{\rho^3\}=1$, and $\rho^4=0$.
\item $\mathcal{M}$ is not Jacobi--Videv.
\item $\alpha_\phi:=\phi^{\prime}\phi^{\prime}\{\phi^{\prime\prime}\}^{-2}$ is a local isometry invariant of $\mathcal{M}$.
\item The following assertions are equivalent:
\begin{enumerate}
\item $\mathcal{M}$ is curvature homogeneous.
\item $\mathcal{M}$ is isometric to $\mathcal{M}_b$ which is defined by $\phi(y)=e^{by}$ for $b>0$.
\item $\mathcal{M}$ is homogeneous.
\end{enumerate}
\end{enumerate}
\end{theorem}

Let $\vec x=(x_1,x_2,x_3,x_4)$ be the canonical coordinates on $\mathbb{R}^4$. 
One says that a pseudo-Riemannian manifold $\mathcal{M}$ of signature $(2,2)$ is a {\it Walker
manifold} if it admits a parallel totally isotropic $2$-plane field -- see \cite{CGM05} for further
details. Such a manifold is locally isometric to an example of the following form:
\begin{equation}\label{eqn-1.a}
\begin{array}{ll}
g(\partial_{x_1},\partial_{x_3})=g(\partial_{x_2},\partial_{x_4})=1,&
  g(\partial_{x_3},\partial_{x_3})=g_{33}(\vec x),\\
 g(\partial_{x_3},\partial_{x_4})=g_{34}(\vec x),&g(\partial_{x_4},\partial_{x_4})=g_{44}(\vec
x)\,.
\end{array}\end{equation}
There are Jacobi--Videv manifolds of signature $(2,2)$ which arise in the context of Walker geometry where we take $g_{33}=g_{44}=0$. We refer to
\cite{BGGV07} for the proof the following result and also for further information concerning Walker geometry:

\begin{theorem}\label{thm-1.4} 
Let $\mathcal{M}:=(\mathbb{R}^4,g)$ be given by Eq.~(\ref{eqn-1.a}) where $g_{33}=g_{44}=0$. Then $\mathcal{M}$ is
Jacobi--Videv if and only if $g_{34}=x_1P(x_3,x_4)+x_2Q(x_3,x_4)+S(x_3,x_4)$ where either 
\begin{enumerate}
\item $P_{/3}=Q_{/4}$, i.e. $Qdx_3+Pdx_4$ is a closed $1$-form, or
\item There exist $(a,b,c)\ne(0,0,0)$ so $P=\frac{c}{a+bx_3+cx_4}$ and $Q=\frac{b}{a+bx_3+cx_4}$.
\end{enumerate}
\end{theorem}

We remark that such a manifold is Einstein if and only if (2) holds. Thus Jacobi--Videv manifolds which
are not Einstein may be created by taking $(P,Q)$ to satisfy (1) but not (2); these will satisfy $\rho$ is nilpotent but $\rho$ need vanish
identically.

We say that $\mathcal{M}$ is {skew--Videv} if $\mathcal{R}(x,y)\rho=\rho\mathcal{R}(x,y)$ for all $x,y$. The following observation, which is of
interest in its own right, will be crucial in our discussion:
\begin{theorem}\label{thm-1.5}
Let $\mathcal{M}$ be a pseudo-Riemannian manifold. The following assertions are equivalent:
\begin{enumerate}
\item $R(\rho\xi_1,\xi_2,\xi_3,\xi_4)=R(\xi_1,\rho\xi_2,\xi_3,\xi_4)=R(\xi_1,\xi_2,\rho\xi_3,\xi_4)=R(\xi_1,\xi_2,\xi_3,\rho\xi_4)$ for all $\xi_i\in
V$.
\item $\mathcal{M}$ is skew--Videv.
\item $\mathcal{M}$ is Jacobi--Videv.
\end{enumerate}
\end{theorem}

We say $\mathcal{M}$ is {\it Jacobi--Tsankov} if $\mathcal{J}(\xi_1)\mathcal{J}(\xi_2)=\mathcal{J}(\xi_2)\mathcal{J}(\xi_1)$ for all $\xi_1,\xi_2\in
V$ and that $\mathcal{M}$ is {\it mixed--Tsankov} if $\mathcal{J}(\xi_1)\mathcal{R}(\xi_2,\xi_3)=\mathcal{R}(\xi_2,\xi_3)\mathcal{J}(\xi_1)$ for all
$\xi_1,\xi_2,\xi_3\in V$. As a scholium to the proof of Theorem \ref{thm-1.5}, we will obtain the following
\begin{theorem}\label{thm-1.6} 
Let $\mathcal{M}$ be a pseudo-Riemannian manifold. Then $\mathcal{M}$ is Jacobi--Tsankov if and only if $\mathcal{M}$ is mixed--Tsankov.
\end{theorem}

The examples we have discussed in Theorems \ref{thm-1.2} and \ref{thm-1.4} involved nipotent Ricci operators. We now discuss examples which are not
Einstein, which are Jacobi--Videv, and where
$\rho$ is not nilpotent. To do this, it is convenient to pass to the algebraic setting.  Let
$V$ be a finite dimensional vector space which is equipped with a non-degenerate inner product of signature $(p,q)$.
Let $A\in\otimes^4V^*$ be a $4$-tensor. One says that
$\mathfrak{M}:=(V,\langle\cdot,\cdot\rangle,A)$ is a model if $A$ has the symmetries of the Riemann curvature tensor:
\begin{eqnarray*}
&&A(v_1,v_2,v_3,v_4)=A(v_3,v_4,v_1,v_2)=-A(v_2,v_1,v_3,v_4),\\
&&A(v_1,v_2,v_3,v_4)+A(v_2,v_3,v_1,v_4)+A(v_3,v_1,v_2,v_4)=0\,.
\end{eqnarray*}
The associated curvature operator $\mathcal{A}$ and bilinear Jacobi operator $\mathcal{J}$ are characterized
by the identities:
\begin{eqnarray*}
&&\langle\mathcal{A}(v_1,v_2)v_3,v_4\rangle=A(v_1,v_2,v_3,v_4),\\
&&\langle\mathcal{J}(v_1,v_2)v_3,v_4)=\textstyle\frac12(A(v_3,v_1,v_2,v_4)+A(v_3,v_2,v_1,v_4))\,;
\end{eqnarray*}
the classical quadratic Jacobi operator being given by $\mathcal{J}(v):=\mathcal{J}(v,v)$.

\begin{definition}\label{defn-1.7}
\rm Let $\mathfrak{M}_0:=(V_0,(\cdot,\cdot),A_0)$ be a Riemannian model.
Let $\{e_i\}$ be an orthonormal basis for $V_0$. Let $V_1=V_0^+\oplus V_0^-$ be two copies of $V_0$ with bases $\{e_i^+,e_i^-\}$.  Let
$\mathfrak{M}_1:=(V_1,\langle\cdot,\cdot\rangle,A_1)$ where
\begin{eqnarray*}
&&\langle e_i^+,e_i^+\rangle=1,\quad\langle e_i^-,e_i^-\rangle=-1,\\
&&A_1(e_i^-,e_j^+,e_k^+,e_l^+)=A_1(e_i^+,e_j^-,e_k^+,e_l^+)=A_1(e_i^+,e_j^+,e_k^-,e_l^+)\\
&&\qquad=A_1(e_i^+,e_j^+,e_k^+,e_l^-)=A_0(e_i,e_j,e_k,e_l),\\
&&A_1(e_i^+,e_j^-,e_k^-,e_l^-)=A_1(e_i^-,e_j^+,e_k^-,e_l^-)=A_1(e_i^-,e_j^-,e_k^+,e_l^-)\\
&&\qquad=A_1(e_i^-,e_j^-,e_k^-,e_l^+)=-A_0(e_i,e_j,e_k,e_l)\,.
\end{eqnarray*}
\end{definition}

The following result may be used to construct examples of Jacobi--Videv models with
$\rho^2=-4s^2\operatorname{id}$:

\begin{theorem}\label{thm-1.8}
Let $\mathfrak{M}_0$ be a Riemannian Einstein model with Einstein constant $s$ and let $\mathfrak{M}_1$ be given by Definition
\ref{defn-1.7}. Then $\mathfrak{M}_1$ is a neutral signature model with $\rho_{\mathfrak{M}_1}^2=-4s^2\operatorname{id}$ which is Jacobi--Videv.
\end{theorem}

There are geometric examples of this phenomena. Again, we specialize the metric of Eq.~(\ref{eqn-1.a}) appropriately:

\begin{theorem}\label{thm-1.9}
Let $(x_1,x_2,x_3,x_4)$ be coordinates on
$\mathbb{R}^4$. Let $\mathcal{M}:=(\mathbb{R}^4,g)$ where
$$\begin{array}{ll}
g(\partial_{x_1},\partial_{x_3})=g(\partial_{x_2},\partial_{x_4})=1,&
  g(\partial_{x_3},\partial_{x_4})=s(x_2^2-x_1^2)/2,\\
g(\partial_{x_3},\partial_{x_3})=sx_1x_2,&
g(\partial_{x_4},\partial_{x_4})=-sx_1x_2\,.\vphantom{\vrule height 12pt}
\end{array}$$
Then $\mathcal{M}$ is locally symmetric (i.e. $\nabla R=0$), $\mathcal{M}$ is Jacobi--Videv,  $\mathcal{M}$ is skew--Videv, and
$\rho^2=-s^2\operatorname{id}$.
\end{theorem}

If $s=1$, then $\rho^2=-\operatorname{id}$ so the Ricci operator $\rho$ defines an almost complex structure on $\mathcal{M}$; note that $\rho$ is
self-adjoint with respect to $g$ and not skew-adjoint with respect to $g$ and thus $\rho$ is not unitary. Furthermore, since
$\nabla R=0$, $\nabla\rho=0$. This manifold has been studied independently in a different context by E. Garc\'{\i}a--R\'{\i}o
\cite{CG07}. They have informed us that both $g$ and the associated Ricci tensor are
irreducible $(2,2)$-metrics sharing the same Levi Civita connection and observed that
this cannot happen either in the Riemannian nor the Lorentzian cases following results of \cite{MT93}.

Here is a brief outline to this paper. In Section \ref{sect-2}, we establish Theorem \ref{thm-1.2}, in Section \ref{sect-3}, we establish Theorem
\ref{thm-1.3}, in Section \ref{sect-4} we establish Theorem \ref{thm-1.5}, in Section \ref{sect-5} we establish Theorem \ref{thm-1.8}, and in Section
\ref{sect-6},  we establish Theorem
\ref{thm-1.9}.

\section{The proof of Theorem \ref{thm-1.2}}\label{sect-2}

Let $\psi_{ij/a}:=\partial_{y_a}\psi_{ij}$ and let
$\psi_{ij/ab}:=\partial_{y_a}\partial_{y_b}\psi_{ij}$. The non-zero Christoffel symbols of the first kind are given by:
\begin{eqnarray*}
&&g(\nabla_{\partial_{x_i}}\partial_{x_j},\partial_{y_a})=\psi_{ij/a},\quad\text{and}\\
&&g(\nabla_{\partial_{x_i}}\partial_{y_a},\partial_{x_j})=g(\nabla_{\partial_{y_a}}\partial_{x_i},\partial_{x_j})=-\psi_{ij/a}\,.
\end{eqnarray*}
Let $C^{ab}$ denote the inverse matrix. We adopt the Einstein convention and sum over repeated indices. The non-zero Christoffel symbols of
the second kind are given by:
\begin{eqnarray*}
&&\nabla_{\partial_{x_i}}\partial_{x_j}=C^{cd}\psi_{ij/c}\partial_{y_d},\quad\text{and}\\
&&\nabla_{\partial_{x_i}}\partial_{y_a}=\nabla_{\partial_{y_a}}\partial_{x_i}=-\delta^{kn}\psi_{ik/a}\partial_{\bar x_n}\,.
\end{eqnarray*}
Clearly $\mathcal{R}(\xi_1,\xi_2)\xi_3=0$ if any $\xi_i\in\operatorname{Span}\{\partial_{\bar x_i}\}$. Furthermore,
\begin{eqnarray*}
&&\mathcal{R}(\partial_{x_i},\partial_{x_j})\partial_{x_k}
   =\nabla_{\partial_{x_i}}\nabla_{\partial_{x_j}}\partial_{x_k}-\nabla_{\partial_{x_j}}\nabla_{\partial_{x_i}}\partial_{x_k}\\
&&\qquad=C^{cd}\delta^{rn}\{-\psi_{jk/c}\psi_{ir/d}+\psi_{ik/c}\psi_{jr/d}\}\partial_{\bar x_n},\\
&&\mathcal{R}(\partial_{x_i},\partial_{x_j})\partial_{y_a}
   =\nabla_{\partial_{x_i}}\nabla_{\partial_{x_j}}\partial_{y_a}-\nabla_{\partial_{x_j}}\nabla_{\partial_{x_i}}\partial_{y_a}=0,\\
&&\mathcal{R}(\partial_{x_i},\partial_{y_a})\partial_{x_j}
   =\nabla_{\partial_{x_i}}\nabla_{\partial_{y_a}}\partial_{x_j}-\nabla_{\partial_{y_a}}\nabla_{\partial_{x_i}}\partial_{x_j}
   =-C^{cd}\psi_{ij/ac}\partial_{y_d},\\
&&\mathcal{R}(\partial_{x_i},\partial_{y_a})\partial_{y_b}=
   \nabla_{\partial_{x_i}}\nabla_{\partial_{y_a}}\partial_{y_b}-\nabla_{\partial_{y_a}}\nabla_{\partial_{x_i}}\partial_{y_b}
   =\delta^{kn}\psi_{ik/ab}\partial_{\bar x_n}\\
&&\mathcal{R}(\partial_{y_a},\partial_{y_b})\partial_{x_i}=
  \nabla_{\partial_{y_a}}\nabla_{\partial_{y_b}}\partial_{x_i}-\nabla_{\partial_{y_b}}\nabla_{\partial_{y_a}}\partial_{x_i}\\
&&\qquad=-\nabla_{\partial_{y_a}}\{\delta^{kn}\psi_{ik/b}\partial_{\bar x_n}\}+
   \nabla_{\partial_{y_b}}\{\delta^{kn}\psi_{ik/a}\partial_{\bar x_n}\}=0,\\
&&\mathcal{R}(\partial_{y_a},\partial_{y_b})\partial_{y_c}=0\,.
\end{eqnarray*}
The polarized Jacobi operator  is given by
$$
\mathcal{J}(\xi_1,\xi_2):\xi_3\rightarrow{\textstyle\frac12}\left\{\mathcal{R}(\xi_3,\xi_1)\xi_2+\mathcal{R}(\xi_3,\xi_2)\xi_1\right\}\,.
$$
The {\it Ricci form} $\rho(\xi_1,\xi_2):=\operatorname{Tr}\{\mathcal{J}(\xi_1,\xi_2)\}$ vanishes if any
$\xi_i\in\operatorname{Span}\{\partial_{\bar x_i}\}$.  Set
$$R_{ijk}{}^n=C^{cd}\delta^{rn}\{-\psi_{jk/c}\psi_{ir/d}+\psi_{ik/c}\psi_{jr/d}\}\,.$$
One has:
$$\begin{array}{ll}
  \mathcal{J}(\partial_{x_i},\partial_{x_j})\partial_{x_k}={\textstyle\frac12}(R_{kij}{}^n+R_{kji}{}^n)\partial_{\bar x_n},&
  \mathcal{J}(\partial_{x_i},\partial_{x_j})\partial_{y_a}=C^{cd}\psi_{ij/ac}\partial_{y_d},\\
  \mathcal{J}(\partial_{x_i},\partial_{y_b})\partial_{x_j}=-{\textstyle\frac12}C^{cd}\psi_{ij/bc}\partial_{y_d},&
  \mathcal{J}(\partial_{x_i},\partial_{y_a})\partial_{y_b}=-{\textstyle\frac12}\delta^{kn}\psi_{ik/ab}\partial_{\bar x_n},
   \vphantom{\vrule height 11pt}\\
  \mathcal{J}(\partial_{y_a},\partial_{y_b})\partial_{x_i}=\delta^{kn}\psi_{ik/ab}\partial_{\bar x_n},&
  \mathcal{J}(\partial_{y_a},\partial_{y_b})\partial_{y_c}=0\,.\vphantom{\vrule height 11pt}
\end{array}$$
Thus the non-zero components of the {\it Ricci form} $\rho(\xi_1,\xi_2):=\operatorname{Tr}\{\mathcal{J}(\xi_1,\xi_2)\}$
may be seen to be:
$$\rho(\partial_{x_i},\partial_{x_j})=C^{ac}\psi_{ij/ac}\,.$$
Raising indices shows that the Ricci operator is given by:
\begin{equation}\label{eqn-2.a}
\rho(\partial_{x_i})=C^{ac}\delta^{jk}\psi_{ij/ac}\partial_{\bar x_k},\quad
   \rho(\partial_{y_a})=0,\quad\rho(\partial_{\bar x_i})=0\,.
\end{equation}
Since 
\begin{eqnarray*}
&&\operatorname{Range}\{\mathcal{J}(\xi_1,\xi_2)\}\subset\operatorname{Span}\{\partial_{y_a},\partial_{\bar x_i}\}
  \subset\ker\{\rho\},\\
&&\operatorname{Range}\{\rho\}\subset\operatorname{Span}\{\partial_{\bar x_i}\}\subset
  \ker\{\mathcal{J}(\xi_1,\xi_2)\},
\end{eqnarray*}
one has that $\rho\mathcal{J}(\xi_1,\xi_2)=\mathcal{J}(\xi_1,\xi_2)\rho=0$. Thus $\mathcal{M}$ is Jacobi--Videv.
Since the Ricci operator is nilpotent, $\mathcal{M}$ is pseudo-Einstein. Equation (\ref{eqn-2.a}) shows $\mathcal{M}$ is Einstein if and only if
$\rho=0$, i.e. if $C^{ab}\psi_{ij/ab}=0$ for all $ij$. \hfill\qedbox

\section{The proof of Theorem \ref{thm-1.3}}\label{sect-3}
Let $(x,y,z,\bar x)$ be coordinates on
$\mathbb{R}^4$. Let
$\phi=\phi(y)$ be a smooth function defined on a connected open subset of
$\mathbb{R}$. We consider the metric
$$g(\partial_{x},\partial_{\bar x})=g(\partial_y,\partial_y)=g(\partial_z,\partial_z)=1,\ 
  g(\partial_{x},\partial_z)=2\phi(y)\,.$$
The Christoffel symbols of the first kind are given by:
\begin{eqnarray*}
&&g(\nabla_{\partial_x}\partial_z,\partial_y)=g(\nabla_{\partial_z}\partial_x,\partial_y)=-\phi^\prime,\\
&&g(\nabla_{\partial_x}\partial_y,\partial_z)=g(\nabla_{\partial_y}\partial_x,\partial_z)=
  g(\nabla_{\partial_y}\partial_z,\partial_x)=g(\nabla_{\partial_z}\partial_y,\partial_x)
=\phi^\prime\,.
\end{eqnarray*}
The non-zero covariant derivatives are therefore given by:
\begin{eqnarray*}
&&\nabla_{\partial_x}\partial_y=\nabla_{\partial_y}\partial_x=\phi^\prime\{\partial_z-2\phi\partial_{\bar x}\},\\
&&\nabla_{\partial_y}\partial_z=\nabla_{\partial_z}\partial_y=\phi^\prime\partial_{\bar x},\\
&&\nabla_{\partial_z}\partial_x=\nabla_{\partial_x}\partial_z=-\phi^\prime\partial_y\,.
\end{eqnarray*}
The action of the curvature operator may therefore be described by:
\begin{eqnarray*}
&&\mathcal{R}(\partial_x,\partial_y)\partial_x=\nabla_{\partial_x}\phi^\prime\{\partial_z-2\phi\partial_{\partial_{\bar x}}\}
         =-\phi^\prime\phi^\prime\partial_y,\\
&&\mathcal{R}(\partial_x,\partial_y)\partial_y=-\nabla_{\partial_y}\phi^\prime\{\partial_z-2\phi\partial_{\bar x}\}
         =-\phi^{\prime\prime}\partial_z+\{2\phi^{\prime\prime}\phi+\phi^\prime\phi^\prime\}\partial_{\bar x},\\
&&\mathcal{R}(\partial_x,\partial_y)\partial_z=\nabla_{\partial_x}\phi^\prime\partial_{\bar x}+\nabla_{\partial_y}\phi^\prime\partial_y
         =\phi^{\prime\prime}\partial_y,\\
&&\mathcal{R}(\partial_x,\partial_z)\partial_x=-\nabla_{\partial_x}\phi^\prime\partial_y
         =-\phi^\prime\phi^\prime\{\partial_z-2\phi\partial_{\bar x}\},\\
&&\mathcal{R}(\partial_x,\partial_z)\partial_y=\nabla_{\partial_x}\phi^\prime\partial_{\bar x}
          -\nabla_{\partial_z}\phi^\prime\{\partial_z-2\phi\partial_{\bar x}\}
         =0,\\
&&\mathcal{R}(\partial_x,\partial_z)\partial_z=\nabla_{\partial_z}\phi^\prime\partial_y
         =\phi^\prime\phi^\prime\partial_{\bar x},\\
&&\mathcal{R}(\partial_y,\partial_z)\partial_x=-\nabla_{\partial_y}\phi^\prime\partial_y
         -\nabla_z\phi^\prime\{\partial_z-2\phi\partial_{\bar x}\}
         =-\phi^{\prime\prime}\partial_y,\\ 
&&\mathcal{R}(\partial_y,\partial_z)\partial_y=\nabla_{\partial_y}\phi^\prime\partial_{\bar x}
          =\phi^{\prime\prime}\partial_{\bar x},\\
&&\mathcal{R}(\partial_y,\partial_z)\partial_z=-\nabla_{\partial_z}\phi^\prime\partial_{\bar x}
          =0\,.
\end{eqnarray*}
Consequently
\begin{eqnarray*}
&&\rho(\partial_x,\partial_z)=-\phi^{\prime\prime},\qquad\rho(\partial_x,\partial_x)=2\phi^\prime\phi^\prime,\\
&&\rho:\partial_x\rightarrow-\phi^{\prime\prime}\partial_z+2(\phi^\prime\phi^\prime+\phi\phi^{\prime\prime})\partial_{\bar x},\quad
  \rho:\partial_z\rightarrow-\phi^{\prime\prime}\partial_{\bar x}\,.
\end{eqnarray*}
Since $\rho$ is nilpotent, $\mathcal{M}$ is pseudo-Einstein. We verify that $\mathcal{M}$ is not Jacobi--Videv by computing:
\begin{eqnarray*}
&&\rho\mathcal{J}(\partial_x,\partial_y)\partial_y
   =-\textstyle\frac12\rho\{-\phi^{\prime\prime}\partial_z+(2\phi^{\prime\prime}\phi+\phi^\prime\phi^\prime)\partial_{\bar x}\}=-\textstyle\frac12
\phi^{\prime\prime}\phi^{\prime\prime}\partial_{\bar x},\\
&&\mathcal{J}(\partial_x,\partial_y)\rho\partial_y=0\,.
\end{eqnarray*}

We now study the curvature tensor. We have
\begin{eqnarray*}
&&R(\partial_x,\partial_y,\partial_y,\partial_x)=\phi^\prime\phi^\prime,\quad
R(\partial_x,\partial_z,\partial_z,\partial_x)=\phi^\prime\phi^\prime,\quad
R(\partial_y,\partial_z,\partial_z,\partial_y)=0,\\
&&R(\partial_y,\partial_x,\partial_x,\partial_z)=0,\quad\phantom{...a}
R(\partial_x,\partial_y,\partial_y,\partial_z)=-\phi^{\prime\prime},\quad
R(\partial_x,\partial_z,\partial_z,\partial_y)=0,
\end{eqnarray*}
We say that a basis $\{X,Y,Z,\bar X\}$ is {\it normalized} if
\begin{equation}\label{eqn-3.a}
\begin{array}{lll}
g(X,\bar X)=1,&g(Y,Y)=1,&g(Z,Z)=1,\\
R(X,Y,Y,X)=0,&
R(X,Z,Z,X)=\star,&
R(Y,Z,Z,Y)=0,\bork\\
R(Y,X,X,Z)=0,&
R(X,Y,Y,Z)=-1,&
R(X,Z,Z,Y)=0\,.\bork
\end{array}\end{equation}
Note that $R(X,Z,Z,X)$ is not specified. To create a normalized basis, we set
\begin{eqnarray*}
&&X:=\varepsilon_1\{\partial_x+\delta_1\partial_z-\textstyle\frac12(\delta_1^2+4\phi\delta_1)\partial_{\bar x}\},\\
&&Y:=\partial_y,\quad Z:=\partial_z-(\delta_1+2\phi)\partial_{\bar x},\quad\bar X:=\varepsilon_1^{-1}\partial_{\bar x}\,.
\end{eqnarray*}
We then have:
\begin{eqnarray*}
&&g(X,\tilde X)=g(Z,Z)=g(Y,Y)=1,\\
&&R(X,Y,Y,X)=\varepsilon_1^2\{\phi^\prime\phi^\prime-2\delta_1\phi^{\prime\prime}\},\\
&&R(Y,Z,Z,Y)=R(Y,X,X,Z)=R(X,Z,Z,Y)=0,\\
&&R(X,Y,Y,Z)=-\varepsilon_1\phi^{\prime\prime},\qquad
  R(X,Z,Z,X)=\varepsilon_1^2\phi^\prime\phi^\prime\,.
\end{eqnarray*}
A normalized basis may then be defined by setting:
$$
\varepsilon_1:=\{\phi^{\prime\prime}\}^{-1}\quad\text{and}\quad
\delta_1=\textstyle\frac12\phi^\prime\phi^\prime\{\phi^{\prime\prime}\}^{-1}\,.
$$

We study the group of symmetries. We note that $\rho:X\rightarrow -Z$ and $\rho:Z\rightarrow -\tilde X$. Let $\{X_1,Y_1,Z_1,\bar X_1\}$ be another
normalized model. Equation (\ref{eqn-3.a}) yields
\begin{eqnarray*}
&&V_1:=\operatorname{Span}_{\xi_i\in\mathbb{R}^4}\{\mathcal{R}(\xi_1,\xi_2)\xi_3\}
  =\operatorname{Span}\{Y,Z,\bar X\}=\operatorname{Span}\{Y_1,Z_1,\bar X_1\},\\
&&V_1^\perp=\operatorname{Span}\{\bar X\}=\operatorname{Span}\{\bar X_1\},\\
&&\ker(\rho)=\operatorname{Span}\{Y,\bar X\}=\operatorname{Span}\{Y_1,\bar X_1\}\,.
\end{eqnarray*}
Consequently, we may express:
$$\begin{array}{rr}
X_1=a_1X+a_2Y+a_3Z+a_4\bar X,&
Y_1=b_1Y+b_2\bar X,\\
Z_1=c_1Y+c_2Z+c_3\bar X,&
\bar X_1=d_1\bar X\,.\bork
\end{array}$$
The conditions on the metric tensor yield $b_1=\pm1$, $c_1=0$, and $c_2=\pm1$. The condition $R(X_1,Y_1,Y_1,X_1)=0$ shows $a_3=0$; as
$R(Y_1,X_1,X_1,Z_1)=0$, one also has $a_2=0$. Thus since $g(X_1,Y_1)=g(X_1,Z_1)=0$ we also have $b_2=c_3=0$. Since $g(X_1,X_1)=0$, we also have
$a_4=0$. The relations $g(X_1,\bar X_1)=1$ and $R(X_1,Y_1,Y_1,Z_1)=-1$ then imply $a_1=c_2$ and $d_1=a_1^{-1}$. Thus
$$X_1=a_1X,\quad Y_1=b_1Y,\quad Z_1=a_1Z,\quad\bar X_1=a_1^{-1}\bar X\quad\text{for}\quad a_1^2=b_1^2=1\,.$$

The calculations performed above show that 
$$\alpha_\phi:=R(X,Z,Z,X)= \{\phi^{\prime\prime}\}^{-2}\phi^\prime\phi^\prime$$
is a local isometry invariant of $\mathcal{M}$; Assertion (2) of Theorem \ref{thm-1.3} follows. 

Assume that $\mathcal{M}$ is curvature homogeneous. By Assertion (2), $\alpha_\phi$ is constant. This
implies that
$\phi=ae^{by}+c$ where $a$, $b$, and $c$ are suitably chosen real constants with $a\ne0$ and $b\ne0$ and consequently
$$ds^2=dx\circ d\bar x+dy\circ dy+dz\circ dz+2\{ae^{by}+c\}dx\circ dz\,.$$
We change variables setting $x=a^{-1}x_1$, $y=\operatorname{sign}(b)y_1$, $z=z_1$, and $\bar x=a\bar x_1-2cz$. We show Assertion (3a) implies
Assertion (3b) by checking:
\begin{eqnarray*}
&&ds^2=a^{-1}dx_1\circ(ad\bar x_1-2cdz)+dy_1\circ dy_1+dz_1\circ dz_1\\&&\qquad+2(ae^{|b|y_1}+c)a^{-1}dx_1\circ dz_1\\
&=&dx_1\circ d\bar x_1+dy_1\circ dy_1+dz_1\circ dz_1+2e^{|b|y_1}dx_1\circ dz_1\,.
\end{eqnarray*}

We therefore suppose $\phi(y)=e^{by}$ for $b>0$. Given $(a_1,a_2,a_3,a_4)\in\mathbb{R}^4$, let $T(x,y,z,\bar
x)=(e^{-2ba_2}x+a_1,y+a_2,z+a_3,e^{2ba_2}\bar x+a_4)$. Then
\begin{eqnarray*}
T^*ds^2&=&e^{-2ba_2}dx\circ e^{2ba_2}d\bar x+dy\circ dy+dz\circ dz\\
   &+&2e^{b(y+a_2)}e^{-2ba_2}dx\circ dz=ds^2\,.
\end{eqnarray*}
Consequently, $T$ is an isometry of $\mathcal{M}$ with $T(0,0,0,0)=(a_1,a_2,a_3,a_4)$. This shows that the group of isometries acts transitively on
$\mathbb{R}^4$ which completes the proof of Theorem \ref{thm-1.3}.
\hfill\qedbox
\section{The proof of Theorem \ref{thm-1.5}}\label{sect-4}

We begin with a brief
technical observation:
\begin{lemma}\label{lem-4.1}
Let $\mathfrak{M}$ be a model and let $T$ be a self-adjoint linear map of $V$. The following conditions are
equivalent:
\begin{enumerate}
\item $A(T\xi_1,\xi_2,\xi_3,\xi_4)=A(\xi_1,T\xi_2,\xi_3,\xi_4)=A(\xi_1,\xi_2,T\xi_3,\xi_4)=
A(\xi_1,\xi_2,\xi_3,T\xi_4)$ for all $\xi_i\in V$.
\item $T\mathcal{A}(\xi_1,\xi_2)=\mathcal{A}(\xi_1,\xi_2)T$ for all $\xi_i\in V$.
\item $T\mathcal{J}(\xi)=\mathcal{J}T(\xi)$ for all $\xi\in V$.
\end{enumerate}\end{lemma}

\begin{proof}
Let $\vec\xi:=(\xi_1,\xi_2,\xi_3,\xi_4)\in V^4$. Let $\{i,j,k,l\}$ denote a permutation of the indices $\{1,2,3,4\}$. Set
$a_{ijkl}=a_{ijkl}(\vec\xi):=A(T\xi_i,\xi_j,\xi_k,\xi_l)$; $a$ need not have the symmetries of an algebraic curvature tensor.
Conditions (1), (2), and (3) are equivalent, respectively, to the following identities:
\begin{eqnarray}
&&a_{ijkl}=-a_{jikl}=a_{klij}=-a_{lkij},\label{eqn-X.a}\\
&&a_{ijkl}=-a_{jikl},\label{eqn-X.b}\\
&&a_{ijkl}+a_{ikjl}=a_{ljki}+a_{lkji}\,.\label{eqn-X.c}
\end{eqnarray}
Clearly Eq.~(\ref{eqn-X.a}) implies both Eqs.~(\ref{eqn-X.b}) and (\ref{eqn-X.c}).

Conversely, suppose Eq.~ (\ref{eqn-X.b}) holds. We may express:
$$\begin{array}{lll}
a_{1234}=\alpha_1,&a_{1342}=\alpha_2,&a_{1423}=-\alpha_1-\alpha_2,\\
a_{2134}=-\alpha_1,&a_{2341}=\alpha_3,&a_{2413}=\alpha_1-\alpha_3,\\
a_{3124}=\alpha_2,&a_{3241}=-\alpha_3,&a_{3412}=-\alpha_2+\alpha_3,\\
a_{4123}=\alpha_1+\alpha_2,&a_{4231}=\alpha_1-\alpha_3,&a_{4312}=\alpha_2-\alpha_3\,.
\end{array}$$
The Bianchi identity then implies $\alpha_3=\alpha_1+\alpha_2$ and Eq.~(\ref{eqn-X.a}) then follows.

Finally, suppose Eq.~(\ref{eqn-X.c}) holds. We use the Bianchi identities to express:
$$\begin{array}{lll}
a_{1234}=\alpha_1,&a_{1342}=\alpha_2,&a_{1423}=-\alpha_1-\alpha_2,\\
a_{2134}=-\beta_1,&a_{2341}=\beta_1+\beta_2,&a_{2413}=-\beta_2,\\
a_{3124}=\gamma_2,&a_{3241}=-\gamma_1-\gamma_2,&a_{3412}=\gamma_1,\\
a_{4123}=\delta_1+\delta_2,&a_{4231}=-\delta_2,&a_{4312}=-\delta_1\,.
\end{array}$$
We use Eq.~(\ref{eqn-X.c}) to derive the 6 identities:
\begin{eqnarray*}
&&\alpha_1+2\alpha_2=a_{1342}+a_{1432}=
  a_{2341}+a_{2431}=\beta_1+2\beta_2,\\
&&-2\alpha_1-\alpha_2=a_{1243}+a_{1423}=a_{3241}+a_{3421}=-2\gamma_1-\gamma_2,\\
&&\alpha_1-\alpha_2=a_{1234}+a_{1324}=a_{4231}+a_{4321}=\delta_1-\delta_2,\\
&&\beta_1-\beta_2=a_{2143}+a_{2413}=a_{3142}+a_{3412}=\gamma_1-\gamma_2,\\
&&-2\beta_1-\beta_2=a_{2134}+a_{2314}=a_{4132}+a_{4312}=-2\delta_1-\delta_2,\\
&&\gamma_1+2\gamma_2=a_{3124}+a_{3214}=a_{4123}+a_{4213}=\delta_1+2\delta_2\,.
\end{eqnarray*}
We set $\beta_i=\alpha_i+\varepsilon_i$, $\gamma_i=\alpha_i+\varrho_i$, and $\delta_i=\alpha_i+\sigma_i$.
We then have:
$$\begin{array}{lll}
0=\varepsilon_1+2\varepsilon_2,&
0=-2\varrho_1-\varrho_2,&
0=\sigma_1-\sigma_2,\\
\varepsilon_1-\varepsilon_2=\varrho_1-\varrho_2,&
-2\varepsilon_1-\varepsilon_2=-2\sigma_1-\sigma_2,&
\varrho_1+2\varrho_2=\sigma_1+2\sigma_2\,.
\end{array}$$
We use the first three equations to see $\varepsilon_1=-2\varepsilon_2$, $\varrho_2=-2\varrho_1$, and
$\sigma_1=\sigma_2$. The final 3 equations then become:
$$-3\varepsilon_2=3\varrho_1,\quad3\varepsilon_2=-3\sigma_1,\quad-3\varrho_1=3\sigma_1\,.
$$
These equations imply $\varepsilon_i=\varrho_i=\sigma_i=0$ and hence $\alpha_i=\beta_i=\gamma_i=\delta_i$
which completes the proof of the Lemma by showing Eq.~(\ref{eqn-X.a}) holds.\end{proof}

Theorem \ref{thm-1.5} now follows from Lemma \ref{lem-4.1} by taking $T=\rho$ and Theorem
\ref{thm-1.6} follows by taking $T=\mathcal{J}(x)$.

\section{The proof of Theorem \ref{thm-1.8}}\label{sect-5}

 We complexify and let $V_{\mathbb{C}}:=V_0\otimes_{\mathbb{R}}\mathbb{C}$. We extend $(\cdot,\cdot)$
and $A_0$ to be complex multi-linear. Let $\{e_i\}$ be an orthonormal basis for $V_0$. Let $\{e_i^+:=e_i,e_i^-:=\sqrt{-1}e_i\}$ be a basis for the
underlying real vector space $V_1:=V\oplus\sqrt{-1}V$. Let $\Re$ and $\Im$ denote the real and imaginary parts of a complex number,
respectively. It is then immediate that
$$\langle\cdot,\cdot\rangle:=\Re\{(\cdot,\cdot)\}\quad\text{and}\quad
  A_1(\cdot,\cdot,\cdot,\cdot)=\Im\{A_0(\cdot,\cdot,\cdot,\cdot)\}\,.$$
Consequently, $A_1$ has the appropriate curvature symmetries and defines an algebraic curvature tensor. We study the Ricci tensor by computing:
\begin{eqnarray*}
\rho_1(e_i^+,e_j^+)&=&\sum_k\left\{A_1(e_i^+,e_k^+,e_k^+,e_j^+)-A_1(e_i^+,e_k^-,e_k^-,e_j^+)\right\}\\
&=&\Im\sum_k\left\{A_0(e_i,e_k,e_k,e_j)+A_0(e_i,e_k,e_k,e_j)\right\}=0,\\
\rho_1(e_i^-,e_j^-)&=&\sum_k\left\{A_1(e_i^-,e_k^+,e_k^+,e_j^-)-A_1(e_i^-,e_k^-,e_k^-,e_j^-)\right\}\\
&=&\Im\sum_k\{-A_0(e_i,e_k,e_k,e_j)-A_0(e_i,e_k,e_k,e_j)\}=0,\\
\rho_1(e_i^+,e_j^-)&=&\sum_k\left\{A_1(e_i^+,e_k^+,e_k^+,e_j^-)-A_1(e_i^+,e_k^-,e_k^-,e_j^-)\right\}\\
&=&\Im\sum_k\left\{\sqrt{-1}A_0(e_i,e_k,e_k,e_j)+A_0(e_i,e_k,e_k,e_j)\right\}\\
&=&2\sum_kA_0(v_i,v_k,v_k,v_j)=2\rho_0(v_i,v_j)=2s\delta_{ij}
\end{eqnarray*}
since $\mathfrak{M}_0$ is Einstein. We show $\rho_1^2=-4s^2\operatorname{id}$ by computing:
$$
\rho_1:e_i^+\rightarrow-2se_i^-\quad\text{and}\quad\rho_1:e_i^-\rightarrow2se_i^+\,.
$$

We can view $\rho_1=-2s\sqrt{-1}$ as a complex linear map of $V_{\mathbb{C}}$. Since we extended $A_0$
to be complex multi-linear, we compute
\begin{eqnarray*}
&&A_1(\rho_1x,y,z,w)=\Im\{A_0(\rho_1x,y,z,w)\}=\Im\{A_0(-2s\sqrt{-1}x,y,z,w)\}\\
&=&\Im\{-2s\sqrt{-1}A_0(x,y,z,w)\}=\Im\{A_0(x,y,z,-2s\sqrt{-1}w)\}\\
&=&\Im\{A_0(x,y,z,\rho_1w)\}=A_1(x,y,z,\rho_1w)\,.
\end{eqnarray*}
Theorem \ref{thm-1.5} now shows this model is Jacobi--Tsankov and Jacobi--Videv.\hfill\qedbox

\section{The proof of Theorem \ref{thm-1.9}}\label{sect-6}

Theorem \ref{thm-1.9}
will follow from Theorem \ref{thm-1.5} and from the following result:

\begin{lemma}\label{lem-6.1}
Let $\mathcal{M}$ be the manifold of Theorem \ref{thm-1.9}. Then
\begin{enumerate}
\item $\mathcal{M}$ is locally symmetric.
\item Let $\kappa=\frac s2$. Then
$R_{1314}=\kappa$, $R_{1323}=-\kappa$, $R_{1424}=\kappa$, and $R_{2324}=-\kappa$.
\item $\rho\partial_{x_1}=-s\partial_{x_2}$, $\rho\partial_{x_2}=s\partial_{x_1}$,
      $\rho\partial_{x_3}=s\partial_{x_4}$, and $\rho\partial_{x_4}=-s\partial_{x_3}$.
\item $R(\rho\xi_1,\xi_2,\xi_3,\rho\xi_4)=-s^2R(\xi_1,\xi_2,\xi_3,\xi_4)$ for all $\xi_i$.
\item $R(\rho\xi_1,\xi_2,\xi_3,\xi_4)=R(\xi_1,\xi_2,\xi_3,\rho\xi_4)$ for all $\xi_i$.
\end{enumerate}
\end{lemma}

\begin{proof} We used a Mathematica package developed by M. Brozos-V\'azquez, J.C. D\'iaz-Ramos, E. Garc\'ia-R\'io and R.
V\'azquez-Lorenzo to establish Assertions (1)-(3). We establish Assertion (4) by computing:
$$\begin{array}{lll}
R_{1314}=\phantom{-}\kappa,&
R_{\rho1,\rho3,1,4}=-s^2R_{2414}=-s^2\kappa,&
R_{\rho1,3,\rho1,4}=\phantom{-}s^2R_{2324}=-s^2\kappa\\&
R_{\rho1,3,1,\rho4}=\phantom{-}s^2R_{2313}=-s^2\kappa,&
R_{1,\rho3,\rho1,4}=-s^2R_{1424}=-s^2\kappa\\&
R_{1,\rho3,1,\rho4}=-s^2R_{1413}=-s^2\kappa,&
R_{1,3,\rho1,\rho4}=\phantom{-}s^2R_{1323}=-s^2\kappa,\\
R_{1323}=-\kappa,&
R_{\rho1,\rho3,2,3}=-s^2R_{2423}=\phantom{-}s^2\kappa,&
R_{\rho1,3,\rho2,3}=-s^2R_{2313}=\phantom{-}s^2\kappa,\\&
R_{\rho1,3,2,\rho3}=-s^2R_{2324}=\phantom{-}s^2\kappa,&
R_{1,\rho3,\rho2,3}=\phantom{-}s^2R_{1413}=\phantom{-}s^2\kappa,\\&
R_{1,\rho3,2,\rho3}=\phantom{-}s^2R_{1424}=\phantom{-}s^2\kappa,&
R_{1,3,\rho2,\rho3}=\phantom{-}s^2R_{1314}=\phantom{-}s^2\kappa,\\
R_{1424}=\phantom{-}\kappa,&
R_{\rho1,\rho4,2,4}=\phantom{-}s^2R_{2324}=-s^2\kappa,&
R_{\rho1,4,\rho2,4}=-s^2R_{2414}=-s^2\kappa,\\&
R_{\rho1,4,2,\rho4}=\phantom{-}s^2R_{2423}=-s^2\kappa,&
R_{1,\rho4,\rho2,4}=-s^2R_{1314}=-s^2\kappa,\\&
R_{1,\rho4,2,\rho4}=\phantom{-}s^2R_{1323}=-s^2\kappa,&
R_{1,4,\rho2,\rho4}=-s^2R_{1413}=-s^2\kappa,\\
R_{2324}=-\kappa,&
R_{\rho2,\rho3,2,4}=\phantom{-}s^2R_{1424}=\phantom{-}s^2\kappa,&
R_{\rho2,3,\rho2,4}=\phantom{-}s^2R_{1314}=\phantom{-}s^2\kappa,\\&
R_{\rho2,3,2,\rho4}=-s^2R_{1323}=\phantom{-}s^2\kappa,&
R_{2,\rho3,\rho2,4}=\phantom{-}s^2R_{2414}=\phantom{-}s^2\kappa,\\&
R_{2,\rho3,2,\rho4}=-s^2R_{2423}=\phantom{-}s^2\kappa,&
R_{2,3,\rho2,\rho4}=-s^2R_{2313}=\phantom{-}s^2\kappa\,.
\end{array}$$
Since $\rho^2=-s^2$, Assertion (5) follows from Assertion (4).
\end{proof}

\section*{Acknowledgments} Research of P. Gilkey partially supported by the
Max Planck Institute in the Mathematical Sciences (Leipzig) and by Project MTM2006-01432 (Spain). 
Research of S. Nik\v cevi\'c partially supported by Project 144032 (Srbija).
It is a pleasure to acknowledge useful discussions with Professor Videv concerning these matters.

\end{document}